\documentclass{article}
\usepackage{graphicx} 

\begin{document}

\title{Definition formula for probabilities of independent events as a consequence of the Insufficient Reason Principle}

\author{Alexander Dukhovny\\
\small{Department of Mathematics, San Francisco State University} \\
\small{San Francisco, CA 94132, USA} \\
\small{dukhovny [at] sfsu.edu}\\ }

\date{\today}          

\maketitle

\begin{abstract}

The standard definition formula for probabilities of independent events is derived as a consequence of the Insufficient Reason Principle expressed as the Maximum Relative Divergence Principle for grading (order-comonotonic) functions on a totally ordered set.

\end{abstract}

\section{Introduction}
\label{sec:Intro}

The standard definition formula $P(A\cap B)=P(A)P(B)$ for probabilities of independent events $A, B$ is an integral part of the foundation of the Probability Theory. The formula itself (or its equivalent forms) has been, so far, either directly postulated or based on the intuitive concepts of the Conditional Probability and Independence of events.

In this article we continue the work started in \cite{P(A|B)} to derive that definition formula from an underlying fundamental Insufficient Reason Principle (IRP). In many probability theory problems IRP has been expressed as the Maximum Shannon Entropy Principle (MEP). In that form, it has been used for choosing the IRP-suggested probability distribution (and/or its missing parameters) as the one that maximizes, under the application-specific constraints, the Shannon Entropy functional (see, e.g., \cite{Shannon}, \cite{Jaynes}). That approach has proved extremely effective in very many cases.

The original Shannon Entropy formula and MEP have been generalized in many ways: Kullback-Leibler (K-L) Relative Entropy and Divergence, Partition Entropy, Kolmogorov-Sinai Entropy, Topological Entropy, Entropy of General Measures and many others - see, e.g., the review in \cite{Review}. 

In this paper, similar to \cite{P(A|B)}, we use the approach started in \cite{Dukh} and continued in \cite{AOME} to generalize Shannon Entropy and associated Kullback-Leibler Relative Divergence of probability distributions to that of the Relative Divergence (RD) of one "grading" (order-comonotonic) function from another on a totally ordered set. (The very term Relative Divergence comes from \cite{K-L}.) 

Accordingly, Maximum Relative Divergence Principle is used in place of MEP as the mathematical expression of IRP as a method to obtain the IRP-suggested grading function $F$ on a totally ordered set $W$ ("chain") where a prior "null" grading function $G$ is available.

\section{Basic Definitions and Properties}\label{Basics}

The initial RD setup (see \cite{AOME}) begins as follows: let $W = \{ w_k, \quad k \in Z$\} be a chain totally ordered by a relation $\prec$. A real-valued function $F$ on $W$ is said to be a grading function (GF) on $W$ if it is order-comonotonic, that is,

$w \prec v \iff F(w) < F(v)$ for all $w, v \in W $.

(For example, the "indexing" function $I: I(w_k) = k$ is a "natural" GF on $W$, used as a "null" GF in the absence of any other information.)

For grading functions $F(w)$ and $G(w)$, the relative divergence of $F$ from $G$ over $W$ is defined as

\begin{equation} \label{RD}
\mathcal{D}(F \Vert G) \vert_W = -\ \sum_{k \in Z}
\ln  \left( \frac{f_k}{g_k}\right) f_k, 
\end{equation}

\noindent where

$f_k=\Delta_k F = F(w_k) - F(w_{k-1})$, $ g_k=\Delta_k G = G(w_k) - G(w_{k-1}), \quad  k \in Z$

\noindent are the increments of, respectively, $F$ and $G$ along the chain $W$. (Where $W$ is infinite absolute convergence of the series is to be guaranteed.)

When $F$ is a probability cumulative distribution function on $W$ and $G = I$, equation (\ref{RD}) reduces to Shannon Entropy of that probability distribution:

\begin{equation} \label{RD-Shannon}
\mathcal{D}(F \Vert I) \vert_W = -\ \sum_{k \in Z} f_k \ln {f_k}   , 
\end{equation}

The Maximum Relative Divergence Principle (MRDP) for grading functions on chains is introduced as a generalization of the Maximum Entropy Principle as follows:

MRDP: Suppose a "null" grading function $G$ is defined on a chain $W$. Among all application-admissible grading functions on $W$ with the same value range, $F$ is said to be IRP-suggested (or "least-presuming") if its Relative Divergence from the given $G$ over $W$ is the highest possible.

\newpage \section{Definition formula for probabilities of independent events as a consequence of MRDP}

Consider a random experiment with the outcome space $U$ where events $A, B, A\cap B, A\cup B$ are defined. Applying The Insufficient Reason Principle to make a decision on the independence of $A$ and $B$ their probabilities $p_1 = P(A)$ and $p_2 = P(B)$ must be specified - but no other information is to be presumed.

Proposition 1. In the absence of any other information on the relation of events $A$ and $B$, using the Maximum Relative Divergence Principle as an expression of the Insufficient Reason Principle leads to the classic definition formula for probabilities of independent events: 

\begin{equation}\label{independence}
P(A\cap B) = P(A)P(B)  
\end{equation}

Proof. Consider the powerset of $U$ with subset inclusion as the ordering relation $\prec$. The following maximal event chains emerge on the powerset of $U$:

$C_1 = \{ \emptyset, A\cap B, A, A\cup B, U \}, \quad C_2 = \{ \emptyset, A\cap B, B, A\cup B, U.\}$

In the absence of any probability assignment on $U$, the only natural grading function defined on both chains is the indexing function $G$ whose values on each of the chains are $\{ 0, 1, 2, 3, 4 \}$.

Assigning probabilities to the defined events, grading functions $F_1$ and $F_2$ on the chains $C_1$ and $C_2$ are assigned their values based on the available information. Some values are predetermined by the probability axioms: 

$F_1(\emptyset) = F_2(\emptyset) = 0, \quad F_1(U) = F_2(U) = 1$. 

When, as it is supposed for application of the Insufficient Reason Principle, the new available information specifies only the probabilities $p_1$ and $p_2$, we assign

$F_!(A) = p_1, \quad F_2(B) = p_2$

\noindent Denoting the unknown $P(A\cap B) = x$, by the probability axioms, we assign 

$F_1(A\cap B) = F_2(A\cap B) = x,$ \quad and

$F_1(A\cup B) = F_2(A\cup B) = P(A\cup B) = p_1 + p_2 - x$,

\noindent completing the definitions of $F_1$ and $F_2$,

As defined, $x \in {[max(0, p_1+p_2 - 1), min(p_1, p_2)]}$. 

Also, the increments of $G$ along both $C_1$ and $C_2$ are $g_i = 1, i= 1,2,3,4$. 

Now we follow (\ref{RD}) to compute 

$d_1(x) = \mathcal{D}(F_1 \Vert G) \vert_{C_1} $ 
and $d_2(x) = \mathcal{D}(F_2 \Vert G) \vert_{C_2}$

- the Relative Divergences of $F_1$ and $F_2$ from $G$ over $C_1$ and $C_2$. 

It turns out that $d_1(x) = d_2(x) = d(x) = $

$-x\ln{x}-(p_1-x)\ln{(p_1-x)} - (p_2-x)\ln{(p_2-x)}$ 
$- (1-p_1-p_2 + x)\ln{(1-p_1-p_2 + x)}.$

Using MRDP, we now find the value of $x$ that maximizes $d(x)$ on its domain - the interval 
$max[0, p_1+p_2 - 1]<  x < min[p_1, p_2]$.  It follows that

$d'(x) = - \ln{x} +\ln{(p_1-x)} + \ln{(p_2 - x)} - \ln{(1-p_1 - p_2 +x)}, $

\noindent and the only root of $d'(x)$ in its domain is $x = p_1 p_2.$ Next,

$d''(x) = -\frac{1}{x} - \frac{1}{(p_1 - x)} - \frac{1}{(p_2 - x)}- \frac{1}{(1-p_1 - p_2 +x)} < 0,$

\noindent so $d(x)$ is concave down on its domain, and, therefore, the only maximum of $d(x)$ there is at $x = p_1p_2$. That is, it follows from MRDP that 

$P(A\cap B) = P(A)P(B)$, thereby concluding the proof of Proposition 1.

\end{document}